\documentclass{article}

\usepackage{dcolumn}% Align table columns on decimal point
\usepackage{bm}% bold math

\usepackage{amsthm}
\usepackage{amssymb}
\usepackage{graphicx}

\usepackage{amstext}
\usepackage{amsmath}
\usepackage{amsfonts}
\usepackage{units}
\usepackage{mathrsfs}

\theoremstyle{theorem}
\newtheorem{theorem}{Theorem}

\theoremstyle{theorem}

\theoremstyle{theorem}

\theoremstyle{lemma}

\newtheorem{remark}{Remark}

\theoremstyle{definition}

\theoremstyle{example}

\newcommand{\p}{\partial}

\newcommand{\dd}{{\rm d}}

\begin{document}


\title{Pseudo-Finsler spaces modeled on a  pseudo-Minkowski space}


%\author{A. Garc\'\i a-Parrado G\'omez-Lobo \thanks{F\'{\i}sica Te\'orica, Universidad del Pa\'\i s Vasco,
%Apartado 644, 48080 Bilbao, Spain.
%\\
%E-mail: alfonso@math.uminho.pt
%\footnote{Current address: Institute of Theoretical Physics,Faculty of Mathematics and Physics, Charles University in Prague, V~Hole\v{s}ovi\v{c}k\'ach~2, 180~00 Praha 8,\\ Czech Republic.}
%}}

\author{ A. Garc\'\i a-Parrado G\'omez-Lobo  \\ F\'{\i}sica Te\'orica, Universidad del Pa\'\i s Vasco,\\
Apartado 644, 48080 Bilbao, Spain.
% \footnote{Current address: Institute of Theoretical Physics,Faculty of Mathematics and Physics,\\ Charles University in Prague, V~Hole\v{s}ovi\v{c}k\'ach~2, 180~00 Praha 8, Czech Republic}
\\ e-mail: alfonso@math.uminho.pt  \\[2ex]
         E. Minguzzi  \\ Dipartimento di Matematica e Informatica ``U. Dini'',\\ Universit\`a
degli Studi di Firenze,\\ Via S. Marta 3,  I-50139 Firenze, Italy. \\ e-mail:  ettore.minguzzi@unifi.it}

%\author{A. Garc\'\i a-Parrado G\'omez-Lobo \thanks{F\'{\i}sica Te\'orica, Universidad del Pa\'\i s Vasco,
%Apartado 644, 48080 Bilbao, Spain.  E-mail: alfonso@math.uminho.pt
%}}

%\author{E. Minguzzi\thanks{
%Dipartimento di Matematica e Informatica ``U. Dini'', Universit\`a
%degli Studi di Firenze, Via S. Marta 3,  I-50139 Firenze, Italy.
% E-mail: ettore.minguzzi@unifi.it }}

% \date{2018}

\maketitle

\begin{abstract}
We adopt a vierbein formalism to study pseudo-Finsler spaces modeled on a pseudo-Minkowski space.
 We show that it is possible to obtain closed expressions for most of the geometric objects of the theory, including Berwald's curvature,
 Landsberg's tensor, Douglas' curvature, non-linear connection and Ricci scalar. These expressions are particularly convenient in computations
 since they factor the dependence on the base and the fiber. As an illustration, we
 study Lorentz-Finsler spaces modeled on the Bogoslovsky  Lorentz-Minkowski  space, and give sufficient conditions which guarantee the Berwald property. We then specialize to a recently proposed Finslerian pp-wave metric. Finally, the paper
 points out that non-trivial Berwald spaces have necessarily  indicatrices which admit some non-trivial linear  group of symmetries.
\end{abstract}

%\begin{keyword}


%Berwald space \sep Minkowski space \sep Ricci scalar.

{Keywords: Finsler geometry, Minkowski space, Berwald space.}

%\end{keyword}

%\title{Pseudo-Finsler spaces modeled on a  pseudo-Minkowski space}

%
%
%\date{}
%
%\maketitle
%
%
%\begin{abstract}
% \noindent We adopt a vierbein formalism to study pseudo-Finsler spaces modeled on a pseudo-Minkowski space.
% We show that it is possible to obtain closed expressions for most tensors in the theory, including Berwald,
% Landsberg, Douglas, non-linear connection, Ricci scalar. These expressions are particularly convenient in computations
% since they factor the dependence on the base and the fiber. Finally, the paper contains some comments on the  geometry of Berwald spaces.
%\end{abstract}

\section{Introduction}

In this  work by pseudo-Minkowski space we mean a pair $(V,L)$ where $V$ is an $m$-dimensional vector space
and $L\colon \bar\omega \to \mathbb{R}$ is a function  defined on the closure of an
open cone\footnote{We do not assume it to be convex.} $\omega\subset V$
with vertex the origin, such that
\begin{itemize}
\item[(a)] $L\in C^4(\omega)\cap C^0(\bar \omega)$,
\item[(b)] $\forall s>0$, and $y\in V$ we have $L(s y)=s^2 L(y)$,
\item[(c)] the Hessian $g:=(\p_\mu \p_\nu {L}) \dd y^\mu \otimes \dd y^\nu$ is non-degenerate on $\omega$.
\end{itemize}
It is understood that the Hessian is calculated using any set of linear coordinates on $V$
(namely a dual basis on $V^*$). More specific theories are possible,
for instance if $g$ is positive definite and $\omega= V\backslash 0$ one speaks of
Minkowski space, while if $g$ has Lorentzian signature and $\omega$ is a sharp convex cone one speaks of Lorentz-Minkowski space.
We warn the reader that already in the positive definite case the concept
of Minkowski space is not homogeneously defined in the literature  (compare \cite{martini01,paiva04}).

If we were concerned with just Minkowski and Lorentz-Minkowski spaces we would probably add the condition
 %We mention that shall not use condition (d) in the derivation our results.
\begin{itemize}
\item[(d)] $L\vert_\omega \ne 0$ and $L\vert_{\p \omega}=0$,
\end{itemize}
in the definition. However, we shall not use it in our derivations (so vector spaces endowed with a non-degenerate bilinear form are included in the definition).

A pseudo-Finsler space is instead a pair $(M,\mathscr{L})$ where $M$ is a $m$-di\-men\-sion\-al manifold and
$\mathscr{L}: \bar \Omega \to \mathbb{R}$, $\Omega \subset TM$, $\pi(\Omega)=M$, with the property
that
if we define $\Omega_x:=\pi^{-1}(x)\cap \Omega$,
$\mathscr{L}_x:=\mathscr{L}\vert_{\bar \Omega_x}$, $\Omega_x\subset T_xM$, then the pair
$(T_x M,\mathscr{L}_x)$
is a pseudo-Minkowski
space with a signature independent of $x$.
One can assume various differentiability conditions on the dependence of
$\mathscr{L}$ on $x$, for our calculations $C^2$ will suffice.
Of course, one could define Finsler and Lorentz-Finsler spaces adding the conditions introduced above for pseudo-Minkowski spaces.

\begin{remark}
The domain of $L$ is really of little importance for the purpose of this work since all our calculations are local over $TM$.
We have chosen a closed cone domain because the causality aspects of these theories can be studied with the help of the theory of differential inclusions.
\end{remark}

This paper is concerned with a special type of pseudo-Finsler space,
namely those for which all pseudo-Minkowskian tangent spaces $(T_xM,\mathscr{L}_x)$
are modeled on the same pseudo-Minkowski space $(V, L)$. This means that locally
we can find a linear isomorphism (sufficiently differentiable in $x \in U\subset M$)
\begin{equation}
\varphi_x \colon T_x M \to V,
\end{equation}
such that $\varphi_x(\Omega_x)=\omega$, and $L(\varphi_x(y))=\mathscr{L}_x(y)$ for every $y\in \bar \Omega_x$.
Parallelizable manifolds admit structures of this type and parallelizability can be relaxed provided $(V,L)$
admits a Lie group $G$  of linear isomorphisms (they preserve $\mathscr{L}$ and its domain). In that case
it is sufficient that $M$ admits a $G$-structure \cite{ichijyo76}. Our considerations will be of local
character so we shall not enter into a more detailed discussion \cite{ichijyo76}.
The idea of Finsler space modeled on the same Minkowski space  has been first introduced and investigated by Ichijyo
\cite{ichijyo76} and further progress has been obtained by Aikou in \cite{aikou00}.
Matsumoto \cite{matsumoto78c}, Izumi \cite{izumi83}, Sakaguchi \cite{sakaguchi83},
and Asanov and Kirnasov \cite{asanov84} had also considered similar spaces referring
to them as {\em 1-form Finsler spaces}. Szilasi and Tam\'assy call
them {\em affine deformations of Minkowski spaces} \cite{szilasi12},
while Libing Huang, and Bartelme\ss \ and Matveev borrow from Zhongmin Shen the terminology
 {\em single colored} or {\em monochromatic} Finsler spaces \cite{huang15,bartelmeb17}.

Roughly speaking, while in pseudo-Finsler spaces the anisotropic features change
from point to point and can even be absent in some open subset of $M$,
in pseudo-Finsler spaces modeled on a pseudo-Minkowski space the anisotropic features
do not depend on the point. Of course, pseudo-Riemannian manifolds provide an example
of pseudo-Finsler space modeled on a pseudo-Minkowski space. Berwald spaces are also
modeled on the same Minkowski space. The linear isometry between different tangent spaces
is provided by the Ehresmann non-linear (actually linear) parallel transport connecting the two points. For Berwald spaces the Ehresmann connection is nothing but the pullback of the Berwald connection under any local section.

However, Berwald spaces are  quite rigid.
In fact, let $G_x$ be the  group of linear isomorphisms of
$(T_xM,\mathscr{L}_x)$
and let $\mathfrak{g}_x$ be its Lie algebra. Let $R_x\colon \Lambda^2 T_xM \to \textrm{End}(T_xM)$ be the curvature of the
Ehresmann connection.
We have (for related ideas cf.\ \cite{szabo81})
\begin{theorem}
Let $(M,\mathscr{L})$ be a connected Berwald pseudo-Finsler space.
Then the Lie algebra generated by $R_x$ is contained in $\mathfrak{g}_x$.
Moreover, if for every $x\in M$,  $\mathfrak{g}_x=0$ (e.g.\ because it is modeled on the same pseudo-Minkowski space with no symmetries),
then the  Berwald pseudo-Finsler space is locally an open subset of  pseudo-Minkowski space, while if additionally
$G_x$ is just the identity for every $x$, then it is globally an  open subset of  pseudo-Minkowski space.
\end{theorem}

\begin{proof}
Since the space is Berwald the transport with the Ehresmann connection  is a linear isomorphism of tangent pseudo-Minkowski spaces (it preserves also the Finsler Lagrangian \cite[Prop.\ 10.1.1]{bao00}\cite[Prop.\ 2.1]{minguzzi15}).
Let $X,Y\in T_xM$. By the Ambrose-Singer theorem $R_x(X,Y)$ is an element of the Lie Algebra of $\textrm{Hol}_x$, the holonomy Lie group of the connection at $x$. Thus $t \mapsto \exp  [t  R_x(X,Y)]$
is a one parameter group of linear isomorphisms of
$(T_xM,\mathscr{L}_x)$, hence $R_x(X,Y)\subset \mathfrak{g}_x$.

If for every $x$, $\mathfrak{g}_x=0$ then $R=0$ and since the Ehresmann connection is torsionless it follows that the space is  locally pseudo-Minkowski  (i.e.\ $M$ is covered by patches $U_i$ such that over $TU_i$ in the induced tangent coordinates   $\mathscr{L}$ has no dependence on $x$), see e.g.\ \cite{minguzzi14c}, hence modeled on the same pseudo-Minkowski space. If for every $x$, $G_x$ is just the identity it is an open subset of pseudo-Minkowski space since the transition map in the intersection of two patches must belong to $G$, the linear group of isomorphisms of the model pseudo-Minkowski space, which is trivial.
\end{proof}

This result shows that Berwald spaces  can be non-trivial only if the model pseudo-Minkowski space has a non-trivial linear group of symmetries, a fact which  limits considerably the family of these spaces. The theorem is quite useful, for instance it immediately implies that  a 2-dimensional Berwald Randers space is locally a Riemannian space or an open subset of a Minkowski space, with no need for calculations, in fact the indicatrix is a translated circle and there is no non-trivial  linear Lie  group of symmetries which send the indicatrix to itself unless there is no translation of the center (the same result could have been obtained using \cite[Th.\ 2.4.1 (ii)]{cheng12}).

In recent years there has been a growing interest in Finsler geometry. Unfortunately, in applications the calculations are often prohibitive so it is natural to look for a smaller family of spaces for which they become amenable.

In this work we shall be  interested in pseudo-Finsler  spaces modeled on a Minkowski space which are not necessarily Berwald.
We prove that most tensors of interest can be explicitly calculated. Namely, we shall show that they can be reduced to a  form factored
into terms depending on the base coordinates $x^\mu$ or on the vertical coordinates $y^a$ of the model pseudo-Mikowski space.
This very fact seems to have passed unnoticed, though it implies a dramatic simplification in the investigation of these spaces.
Also it must be remarked that this family, though smaller than the general family of pseudo-Finsler spaces, is in the end fairly large. Most proposed metrics in applications are of this form, we mention the Berwald-Mo{\'o}r metric \cite{berwald39,moor54}, Bogoslovsky and Goenner's metric \cite{bogoslovsky98,bogoslovsky99}, Bogoslovsky's metric \cite{bogoslovsky77,bogoslovsky94}, the conic metric introduced by the second author in \cite{minguzzi16c},  the Randers spaces in which the $\beta$ 1-form has constant $\alpha$-module, or the pp-wave Finsler metric by Fuster and Pabst \cite{fuster16} which we shall more closely study in the last section.

One Finslerian quantity which  attracted our interest is the Ricci scalar. In the theory of Finsler gravity most proposals 
for the  vacuum dynamical equations, starting from Horvath's \cite{horvath50},  imply the condition
\[
\textrm{Ric}(y)=0 ,
\]
which is also suggested by studies of the Raychaudhuri equation \cite{rutz93,minguzzi15}.
Although some of the final tensorial expressions which we shall present are rather concise,
occasionally we had to pass through computations involving hundreds of terms.
This happened precisely in the calculation of $\textrm{Ric}(y)$ so we are happy to have been able
to reduce its expression to a much simpler form, cf.\ Eq. (\ref{ric}).
Here we have been helped by the {\em xAct}
suite of tensor computer algebra \cite{xact}.

\section{The strategy}

Let $\{\tilde e_a\}$ be a basis of $V$ and let $y^a$ be the induced coordinates, so if $\tilde y\in V$,
then there are $y^a \in \mathbb{R}$, $a=1,\ldots,m$, such that $\tilde y=y^a \tilde e_a$
(we adopt the Einstein summation convention).
Let us define at $x\in M$,
\[
y:=\varphi_x^{-1}(\tilde y), \qquad e_a:=\varphi_x^{-1}(\tilde e_a)\;,
\]
so that $y=y^a e_a$.

Let us denote with $\p_\mu$ the basis of $T_xM$ induced by
a local coordinate system $\{x^\mu\}$ on $U\subset M$, then there are invertible matrices
(vierbein) $e^\mu_a$ and $e^a_\mu$, such that $e_a=e_a^\mu \p_\mu$, $e^a=e^a_\mu \dd x^\mu$, where $\{e^a\}$ is the dual basis to $\{e_a\}$,
\begin{equation} \label{map}
e^a_\mu e^\mu_b=\delta^a_b, \qquad e^\mu_a e^a_\nu=\delta^\mu_\nu.
\end{equation}
The vector $y\in T_xM$ can be expanded as follows
\[
y=y^\mu \p_\mu= y^a e_a,
 \]
 where (here we use the physicist convention according to which a different type of index,
 say Greek, Roman,
 might distinguish different objects)
\[
y^\mu=y^a e^\mu_a.
\]
We shall call $y^a$ the {\em vertically induced internal variables}.
 %and $y^\alpha$ the {\em horizontally induced internal variables}.
The commutation coefficients $[e_a,e_b]=c^c{}_{ab} e_c$ are
\[
c^c_{\ ab}= (e^\alpha_a e_{b,\alpha}^\mu-e^\alpha_b e_{a,\alpha}^\mu)\p_\mu=e^\alpha_a e^\beta_b(e^c_{\alpha,\beta}- e^c_{\beta,\alpha}) .
\]
where in the second equality we used (\ref{map}). They can also be obtained through the identity
\begin{equation} \label{nxj}
\dd e^c=\dd (e^c_\gamma \dd x^\gamma)=-c^c_{ab}\, \tfrac{1}{2} e^a \wedge e^b.
\end{equation}
 We shall also denote $c^c_{\ ab,d}=c^c_{\ ab,\mu }e^\mu_d$.

Now according to the isomorphism assumption  between $(T_xM,\mathscr{L}_x)$ and $(V, L)$ we have
\[
\mathscr{L}_x(y^a e_a)=L(y^a \tilde e_a) .
\]
%To simplify notation we work in the next formulas as if $\mathscr{L}_x$ and $L$ depend only on $\{y^a\}$.
In what follows $\tilde e_a$ is considered chosen and
fixed once and for all, so we can regard $L$ as a function
$L\colon\mathbb{R}^m\to \mathbb{R}$. So we denote
\begin{equation} \label{doa}
y_a=\p_a L, \quad g_{ab}=\p_a \p_b L,\quad  C_{abc}=\tfrac{1}{2}\p_a \p_b \p_c L, \quad I_a=g^{bc} C_{abc},  \quad C_{abcd}=\p_aC_{bcd},
\end{equation}
called respectively, the Legendre dual of $y$, the (pseudo-)Finsler metric,  the Cartan torsion, the mean Cartan torsion
and the Cartan curvature of the tangent pseudo-Minkowski space. Since the metric is
non-degenerate we shall also have the inverse metric $g^{ab}$.
To simplify notations we shall also regard $\mathscr{L}$ as a local function
in the horizontally induced variables i.e.\ from $U\times \mathbb{R}^m \to \mathbb{R}$,
$U\subset \mathbb{R}^m$,
so that $(x^\mu,y^\nu)\to \mathscr{L}(x^\mu,y^\nu)$. Thus we can write
\begin{align*}
\mathscr{L}(x^\mu,y^\nu)&= L(e^a_\nu(x) y^\nu).
\end{align*}
Now, most references of Finsler geometry provide expressions of the most important tensors
in the coordinates $(x^\mu, y^\nu)$ induced over the tangent bundle by the coordinates on the base,  so we shall stick to those \cite{minguzzi14c,bejancu00,shen01,szilasi14}.
However, our aim is to change variables from the pair $(x^\mu, y^\nu)$ to $(x^\mu, y^a)$.
So we shall need the change of partial derivatives
\begin{align}
\left.\frac{\p }{\p y^\alpha} \right)_{x^\gamma}&=e^a_\alpha \left.\frac{\p }{\p y^a} \right)_{x^\gamma} , \label{one} \\
\left.\frac{\p }{\p x^\alpha} \right)_{y^\gamma}&=\left.\frac{\p }{\p x^\alpha} \right)_{y^c} +e^b_{\beta,\alpha} e^\beta_c y^c \left.\frac{\p }{\p y^b} \right)_{x^\gamma} . \label{two}
\end{align}
The idea is to arrange the notable tensors in terms of the quantities
$e^a_\alpha$, $e^\alpha_a$, $c^a_{\ bc}$, $c^a_{\ bc,d}$ which only depend on $x^\mu$
and in terms of the quantities displayed in (\ref{doa}) including the inverse metric $g^{ab}$
which only depend on $y^a$. So we do not raise or lower indices with the metric in the commutation
coefficients otherwise we would spoil this property introducing objects dependent on both variables.
Only after we have obtained all expressions of interest we simplify them lowering or raising the indices
of the commutators, with the convention that the upper index in  $c^a_{\ bc}$ gets lowered to the left as $c_{abc}$.
For instance applying the previous formulas to $\mathscr{L}$ we get
\begin{align*}
y_\alpha:=\frac{\p \mathscr{L}}{\p y^\alpha} &= y_a e^a_\alpha, \\
\frac{\p \mathscr{L}}{\p x^\alpha} &= y_b  e^b_{\beta,\alpha} e^\beta_c y^c.
\end{align*}
These derivatives enter the calculation of the spray \cite[Eq.\ (5.2)]{shen01}
\begin{align}
G^\mu&:=\tfrac{1}{2} g^{\mu \sigma} \left(\frac{\p^2 \mathscr{L}}{\p x^\rho \p y^\sigma} y^\rho-\frac{\p \mathscr{L}}{\p x^\sigma}\right) \nonumber\\
%&=\frac{1}{2} g^{\mu \sigma}\left(y_m e^m_{\sigma,\rho} e^\rho_n y^n+g_{cm} e^m_\sigma e^c_{\delta,\rho} e^\delta_d e^\rho_p y^d y^p-y_n e^n_{\beta,\sigma} e^\beta_r y^r\right)\\
&=e^\mu_m \left[\tfrac{1}{2}  \left( g^{m k} y_q y^n c^q_{\ kn}+ e^m_{\delta,\rho} e^\delta_d e^\rho_p y^d y^p\right) \right] . \label{spr}
\end{align}
Observe that  all terms either depend on $x^\mu$ or $y^a$; this fact makes the next calculations easier by repeatedly using (\ref{one})-(\ref{two}).
%However, we give the final expression in the shorter  form
%\begin{equation}
%G^\mu=e^\mu_m \left[\tfrac{1}{2}  \left(  y_q y^n c^{qm}{}_{n}+ e^m_{\delta,\rho} e^\delta_d e^\rho_p y^d y^p\right)\right] ,
%\end{equation}
%which, we stress once again, it is not our working expression for the next %calculations.

The coefficients of the non-linear connection  are
\begin{align}
\begin{split}
N^\mu_\alpha &:=\frac{\p G^\mu}{\p y^\alpha} = e^\mu_l e^b_\alpha
\Big[e^l_{\rho,\eta} e^\eta_b  e^\rho_p y^p +c^m_{\ kn} \left( - C^{l k}_b   y_m y^n  \right. \\
&\qquad \qquad  \qquad \quad \ \left. + \tfrac{1}{2}   \left( g^{l k} g_{mb} y^n+ g^{l k} y_m \delta^n_b+
\delta^l_m \delta^k_b  y^n\right) \right) \Big] .
\end{split}
\end{align}
The coefficients of the Berwald Finsler connection are
\begin{align}
2 G^\rho_{ \alpha \beta}&:=2 \frac{\p^2 G^\rho}{\p y^\alpha \p y^\beta} =[e^\rho_b(e^b_{\alpha, \beta}+e^b_{\beta,\alpha})] \nonumber \\
&\quad + c^e_{\ m n} \Big( g^{r m} g_{e a} \delta^n_b+ g^{r m} g_{e b} \delta^n_a +2 g^{r m} C_{e a b} y^n -2 C^{r m}_{b} g_{e a} y^n \\
&\quad -2 C^{r m}_a g_{e b} y^n-2 C^{r m}_a y_{e}  \delta^n_b   -2 C^{r m}_b y_{e} \delta^n_a -2 \big(\tfrac{\p}{\p y^b} C^{r m}_a\big) y_{e} y^n\Big) . \nonumber
\end{align}
We could simplify the vertical derivative with
\begin{equation}
-\frac{1}{2}\tfrac{\p}{\p y^b} \tfrac{\p}{\p y^a}  g^{r m}=\tfrac{\p}{\p y^b} C^{r m}_a=C^{r m}_{a b}-2 C^{r s }_{b} C_{s a}^m-2 C^{m s }_{b} C^r_{s a} ,
\end{equation}
but the expression gets longer.
The Berwald curvature tensor is defined through  $G^\rho_{ \alpha \beta \gamma}=\frac{\p^3 G^\rho}{\p y^\alpha \p y^\beta \p y^\gamma}$ and in our case it is given by
\begin{align}
\begin{split}
2 G^r_{ a b c}&=2 G^\rho_{ \alpha \beta \gamma} e^r_\rho e^\alpha_a e^\beta_b e^\gamma_c=c^e_{\ m n} \tfrac{\p^3}{\p y^a\p y^b\p y^c} (g^{rm} y_e y^n)  . \label{ber}
%\\
%&=c^e_{\ m n} \Big( 2 g^{r m} C_{e b c} \delta^n_a  +2g^{r m} C_{e c a} \delta^n_b+2 g^{r m} C_{e a b} \delta_c^n\\
%& \qquad \qquad   -2 C^{r m}_a g_{e c}  \delta^n_b -2 C^{r m}_a g_{e b} \delta_c^n -2 C^{r m}_{b} g_{e a} \delta_c^n \\
%& \qquad \qquad -2 C^{r m}_b g_{e c}  \delta^n_a  -2 C_c^{r m} g_{e b} \delta^n_a-2 C_c^{r m} g_{e a} \delta^n_b\\
%&\qquad \qquad -4 C^{r m}_c C_{e a b} y^n-4 C^{r m}_{b} C_{e a c} y^n-4 C^{r m}_a C_{e b c} y^n \\
% &\qquad \qquad  -2 \big(\tfrac{\p}{\p y^c}C^{r m}_{b}\big) g_{e a} y^n  -2 \big(\tfrac{\p}{\p y^c}C^{r m}_a\big) g_{e b} y^n -2 \big(\tfrac{\p}{\p y^b} C^{r m}_a\big) g_{e c} y^n \\
%&\qquad \qquad  -2 \big(\tfrac{\p }{\p y^b}C^{r m}_c\big) y_{e} \delta^n_a-2  \big(\tfrac{\p}{\p y^a}C^{r m}_c\big) y_{e }  \delta^n_b-2 \big(\tfrac{\p}{\p y^a} C^{r m}_b\big) y_{e}  \delta_c^n \\
%&\qquad \qquad +2 g^{r m} C_{e a b c} y^n  + \big(\tfrac{\p^3}{\p y^a\p y^b\p y^c} g^{r m} \big) y_{e}y^n \Big) .
\end{split}
\end{align}
Here it is not convenient to expand the derivatives.
The mean Berwald curvature $E_{bc}:=\frac{1}{2} G^r_{ r b c}$ follows from
\begin{align}
 G^r_{ r b c}&=- c^e_{\ m n}  \tfrac{\p^2}{\p y^b \p y^c} ( I^m y_ e y^n)  . \label{mer}
 &
\end{align}
The Douglas curvature is defined by
$D^\rho_{ \alpha \beta \gamma}=\frac{\p^3 }{\p y^\alpha \p y^\beta \p y^\gamma}
(G^\rho\!-\!\tfrac{1}{m+1} N^a_a y^\rho)$ and in our case it is given by
\begin{align}
 D^r_{ a b c}
 %&= G^r_{ a b c}+\tfrac{c^e_{mn}}{m+1}  \tfrac{\p^3 }{\p y^a \p y^b \p y^c} (I^m y^n y_e y^r) ,\\
 &=c^e_{\ pq} \tfrac{\p^3}{\p y^a\p y^b\p y^c} ((\tfrac{1}{2}g^{rp} +\tfrac{1}{m+1} I^p y^r)y_e y^q)
\end{align}
where $m$ is the dimension of the Finsler space.
In the next expression  $\nabla^{VC}$ is the vertical Cartan derivative
(whose connection coefficients in coordinates $\{y^a\}$ are $C^c_{a b}$) and
$h^a_b:= \delta^a_b-\tfrac{1}{g_y(y,y)} y^a y_b$, $g_y(y,y):=2L$
is the usual projection on the space tangent to the indicatrix.

The Landsberg tensor is
\begin{align}
\begin{split}
L_{a b c}&=L_{\alpha \beta \gamma}  e^\alpha_a e^\beta_b e^\gamma_c=-\tfrac{1}{2} y_r G^r_{ a b c}=-\tfrac{1}{4}c^e_{\ m n} y_r \tfrac{\p^3}{\p y^a\p y^b\p y^c} (g^{rm} y_e y^n) \\
%&=-\frac{1}{2} c^e_{\ m n} \Big( y^{ m} C_{e b c} \delta^n_a  +y^{ m} C_{e c a} \delta^n_b+ y^{ m} C_{e a b} \delta_c^n \\
% &\quad  + C^{m}_{  b c} g_{e a} y^n  + C^{ m}_{c a} g_{e b} y^n+  C^{ m}_{a b}  g_{e c} y^n \\
%&\quad + C^{ m}_{b c} y_{e } \delta^n_a +  C^{ m}_{c a} y_{e } \delta^n_b+ C^{ m}_{a b} y_{e } \delta_c^n   +  \tfrac{1}{2} y_\rho\big(\tfrac{\p^3}{\p y^a\p y^b\p y^c} g^{\rho m} \big) y_{e}  y^n \Big) \\
&=-\tfrac{1}{2} c^e_{\ m n} \Big( y^{ m} C_{e b c} \delta^n_a  +y^{ m} C_{e c a} \delta^n_b+ y^{ m} C_{e a b} \delta_c^n \\
 &\quad  + C^{m}_{  b c} g_{e a} y^n  + C^{ m}_{c a} g_{e b} y^n+  C^{ m}_{a b}  g_{e c} y^n \\
&\quad + C^{ m}_{b c} y_{e } \delta^n_a +  C^{ m}_{c a} y_{e } \delta^n_b+ C^{ m}_{a b} y_{e } \delta_c^n  \\
&\quad   +  2 \big(C^m_{a b c}- C^{m}_{ s a} C^s_{ b c}-  C^m_{s c} C^s_{a b} - C^{m}_{ s b} C^s_{  c a} \big) y_{e}  y^n \Big) \\
&=k_{m  n e} y^{ e} ( C^m_{ b c} \delta^n_a + C^m_{ c a} \delta^n_b+  C^m_{ a b} \delta_c^n) \\
&\quad +  k_{m n e}  y^n y^e  \big(C^m_{ a b c}- C^s_{  b c} C^m_{a s}- C^s_{c a} C^{m}_{b s}  -   C^s_{a b } C^m_{ c s}\big) \label{ssr} \\
&=k_{m  n e} y^{ e} \left[  C^{[m}_{ b c} h^{n]}_a + C^{[m}_{ c a} h^{n]}_b+  C^{[m}_{ a b} h_c^{n]}+ y^{[n} h^{m] s }
h^{p}_a h^{q}_b h^{r}_c \nabla^{VC}_{s} C_{p q r} \right] \; .
\end{split}
\end{align}
The first expressions are particularly suited for calculations since the whole $x$ dependence is contained
in the first  $c^e_{\ m n} $ factor, while the subsequent factors depend only on $y^a$. In the
last two expressions we introduced the object antisymmetric in the first two indices (they are a sort of Ricci rotation coefficients)
\[
 k_{m n e}(x,y)=\tfrac{1}{2}(c_{m  n e}+c_{n   e m} - c_{e m n}) .
\]
The last expression in Eq.\ (\ref{ssr}) is useful because through contraction with $y_n$
it shows that for $m\ge 3$ the tensor in square bracket vanishes if and only if $C_{abc}=0$.
Thus it is not possible to find vertical conditions weaker than $C_{abc}=0$
(i.e.\ the pseudo-Riemannian case)
which guarantee that $L_{abc}=0$.

Still the expression for $L$ simplifies considerably when
$h^{p}_a h^{q}_b h^{r}_c  h^{s}_d \nabla^{VC}_{s} C_{pqr} =0$
which happens if and only if the geometry of the indicatrix is such that $\nabla^{\mathfrak{h}}\mathfrak{c}=0$
where $\mathfrak{c}$ is the Pick cubic form and $\mathfrak{h}$ is the affine metric for the centroaffine transverse
(for translation between Finsler  and centroaffine geometry
the reader is referred to \cite{minguzzi15e}). For instance, $\nabla^{\mathfrak{h}}\mathfrak{c}=0$ holds true whenever
the indicatrix is a homogeneous affine sphere \cite{hu11,hildebrand15}.
%In particular, it is known \cite{hu11,hildebrand15} that $ \nabla^h c=0$ for homogeneous affine spheres for which the previous expression gets considerably simplified.

The mean Landsberg curvature is
\begin{align}
%\begin{split}
J_c:=L^a_{ac}&=-\frac{1}{2} c^e_{m n} \Big( y^mI_e \delta^n_c+I^m g_{ec} y^n+I^m y_e \delta^n_c+2 \left(\tfrac{\p }{\p y^c} I^m +C^m_{sc} I^s
\right)
y_e y^n \Big) \nonumber \\
&=k_{m  n e} y^{ e} \Big[ I^m \delta_c^n  +   y^n \left(\tfrac{\p }{\p y^c} I^m +C^m_{sc} I^s\right)\Big]\label{hap}\nonumber \\
&=k_{m  n e} y^{ e} \Big[ h^{[n}_c I^{m]}+y^{[n} (h^a_c h^{m]}_b \nabla^{VC}_a I^b)\Big] .
%\end{split}
\end{align}
Again, the first expression is useful in calculations since it factors the $x^\mu$ and $y^a$ dependence.

The last formula can be used to show that the vertical tensor in square brackets vanishes if and only if $I=0$ (for the `only if' direction contract first with $y_n$).
Thus it is not possible to find vertical conditions weaker than $I_c=0$ which guarantee that $J_c=0$
(a fact known to hold in general for any pseudo-Finsler space).

The curvature of the non-linear connection is
\[
R^\mu_{\alpha \beta}:= \frac{\delta N^\mu_\beta}{\delta x^\alpha}-\frac{\delta N^\mu_\alpha}{\delta x^\beta}, \qquad \frac{\delta}{\delta x^\alpha}:=\frac{\p}{\p x^\alpha}-N^\mu_\alpha \frac{\p}{\p y^\mu}.
\]
Our calculation for the Ricci scalar gives %\cite{shen01}
\begin{align}
\begin{split}
\textrm{Ric}(y):=& \, R^\mu_{ \mu \alpha} y^\alpha= 2 \frac{\p G^\mu}{\p x^\mu}-y^\mu \frac{\p^2 G^\nu}{\p x^\mu \p y^\nu}+2 G^\mu \frac{\p^2 G^\nu}{\p y^\mu \p y^\nu}-\frac{\p G^\mu}{\p y^\nu}\frac{\p G^\nu}{\p y^\mu}\\
=& \, g^{ad} y^{b} y_{c}
\mathit{c}^{c}{}_{ab,d} + y^{a} y^{b} \mathit{c}^{c}{}_{bc,a}
+I^{b} y^{a} y^{c} y_{d}\mathit{c}^{d}{}_{bc,a} \\
 & + \mathit{c}^{i}{}_{jk}
\mathit{c}^{l}{}_{mn} g_{ia} g_{lb}\Big( \!\tfrac{1}{4}  g^{jm} g^{kn}
y^{a} y^{b} +  g^{b n} g^{km} y^{a} y^{j} - \tfrac{1}{2}
g^{a n} g^{b k} y^{j} y^{m} \\ & - \tfrac{1}{2} g^{ab}
g^{kn} y^{j} y^{m}  + 2 y^{a} y^{j} y^{m}
C^{b k n} - y^a y^{j} y^{b} y^{m}
C^{kg}_h C^{n h}_{g}
\\ & +  I^{j}
 g^{kn} y^{a} y^{b} y^{m}  -  g^{b k} I^{n} y^{a} y^{j} y^{m}   -  g^{en} y^{a} y^{j}
y^{b} y^{m} \tfrac{\partial }{\partial y^e} I^{k}\Big) . \label{ric}
\end{split}
\end{align}
It can be easily checked that if the metric does not depend on the vertical variable ($C_{abc}=0$)
we get the usual expression of the Ricci tensor in terms of the commutation coefficients \cite[Sect.\ 98]{landau62}.

\section{The Bogoslovsky model Lorentz-Minkowski space}
A Lorentz-Minkowski space that has attracted considerable interest is Bogoslovsky's \cite{bogoslovsky94}
 namely $L=-\frac{1}{2} F^2$, where
\begin{equation} \label{mdp}
F=\Big(2 y^0 y^n-\sum^{n-1}_{J=1}(y^J)^2\Big)^{(1-b)/2}(y^0)^b \;,
\end{equation}
and $b$ is a dimensionless parameter. Its study for $n=3$
was revived by the Very Special Relativity theory by Cohen and Glashow  \cite{cohen06}, who observed that most observations are really compatible with a subgroup $SIM(2)$ of the Lorentz group, and by Gibbons, Gomis and Pope \cite{gibbons07} who showed that among the possible continuous deformations $DISIM_b(2)$ is particularly interesting and leads to the local Lagrangian of a Lorentz-Finsler geometry based on the Bogoslovsky Lorentz-Minkowski space.

In this section we apply our method to the Bogoslovsky metric, and in a next section we specialize to a recent pp-wave Lorentz-Finsler space based on it.
It will be convenient to set
\begin{equation} \label{aaa}
A=2 y^0 y^n-\sum^{n-1}_{J=1}(y^J)^2=2 y^0 y^n-\delta_{JK}y^J y^K\;,
\end{equation}
In what is to follow capital Latin indices $J,K,L$ run from $1$ to $n-1$ and $\delta_{JK}$ (resp.\ $\delta^{JK}$) represents the
Euclidean metric (resp.\ its inverse) in dimension $n-1$.
We proceed now with the vertical derivatives.
The differential of $L=-F^2/2$ is
\begin{align}
\begin{split}
y_{0} =& \tfrac{F^{2}}{Ay^0} (-bA+(b-1) y^0y^n)\;,\\
y_J =& \tfrac{F^{2}}{A} (1-b) y^J\;,\\
y_n =& -\tfrac{F^{2}}{A} (1-b) y^0\;.
\end{split} \label{yaa}
 \end{align}
The components of the metric are
\begin{align}
\begin{split}
g_{00}&=b  \left(\tfrac{F}{A y^0}\right)^2[2(1-b) (y^0y^n)^2+4(b-1)y^0y^nA+(1-2b) A^2] \;,\\
g_{0J}&=2(b-1) b \tfrac{F^2}{A^2 }  \tfrac{y^J}{y^0} \left(  y^0 y^n-A\right)\;,\\
g_{0n}&=(1 - b) \tfrac{F^2}{A^2 } \Bigl(-(1 + 2 b)A + 2 b y^0y^n\Bigr)\;,\\
g_{JK}&=(1-b) \tfrac{F^2}{A^2 }  \Bigl(2 b y^J y^K +A\delta_{JK}\Bigr)\;,\\
%g_{JK}&={2 (1 - b) b \tfrac{F^2}{A^2 }  y^J y^K}\;,\\
g_{Jn}&=2 (b-1) b \tfrac{F^2}{A^2 }  y^0 y^J\;,\\
g_{nn}&=2 (1 - b) b \tfrac{F^2}{A^2 } (y^0)^2 .
\end{split} %\nonumber
\end{align}
The components of its inverse are
\begin{align}
\begin{split}
g^{00}&=-\tfrac{2 b (y^0)^2}{(1 + b) F^2}\;,\qquad\\
g^{0J}&=-\tfrac{2 b y^0 y^J}{(1 + b) F^2}\;,\qquad\\
g^{0n}&=-\tfrac{A}{(1 + b) F^2}-\tfrac{ 2 b y^0 y^n}{(1 + b) F^2}\;,\\
g^{JK}&=
\tfrac{A}{(1-b) F^2} \delta^{JK}-\tfrac{2by^J y^K}{(1+b)F^2}\;,\\
%g^{JK}&=-\tfrac{2 b y^J y^K}{(1 + b) F^2}\;,\\
g^{Jn}&=\tfrac{2 b y^J A}{(1-b^2) F^2 y^0}-\tfrac{ 2 b y^J y^n}{(1 + b) F^2} \;,\\
g^{nn}&=\tfrac{b A (A-4y^0y^n) }{(b^2-1) F^2 (y^0)^2} -\tfrac{ 2 b  (y^n)^2}{(1 + b) F^2}  .
\end{split} \label{sof}
\end{align}
%
%\begin{align*}
%&g^{00}g_{00}+g^{0J} g_{J0}+g^{0n}g_{n0}\\
%&=-\tfrac{2 b (y^0)^2}{(1 + b) F^2} b  \tfrac{F^2}{A^2 }  (y^0)^{-2}   [2(1-b) (y^0y^n)^2+4(b-1)y^0y^nA+(1-2b) A^2]\\
%&-\tfrac{2 b y^0 y^J}{(1 + b) F^2}2(b-1) b \tfrac{F^2}{A^2 }  \tfrac{y^J}{y^0} \left(  y^0 y^n-A\right)+[-\tfrac{A}{(1 + b) F^2}-\tfrac{ 2 b y^0 y^n}{(1 + b) F^2}](1 - b) \tfrac{F^2}{A^2 } \Bigl(-(1 + 2 b)A + 2 b y^0y^n\Bigr)\\
%&=-\tfrac{2 b}{(1 + b) } b  \tfrac{1}{A^2 }    [2(1-b) (y^0y^n)^2+4(b-1)y^0y^nA+(1-2b) A^2]\\
%&-\tfrac{2 b y^0 y^J}{(1 + b) }2(b-1) b \tfrac{1}{A^2 }  \tfrac{y^J}{y^0} \left(  y^0 y^n-A\right)+[-\tfrac{A}{(1 + b) F^2}-\tfrac{ 2 b y^0 y^n}{(1 + b) F^2}](1 - b) \tfrac{F^2}{A^2 } \Bigl(-(1 + 2 b)A + 2 b y^0y^n\Bigr)
%\end{align*}
\begin{remark}
The components of the Cartan torsion are too many to be displayed here.
Here we just give the
mean Cartan torsion for $n=3$
\begin{align}
\begin{split}
I_0&={4 b }{(y^0)^{-1} }\;, \quad
I_1={4 b y^1}/A \;, \quad
I_2={4 b y^2}/A \;,\quad
I_3= {-4 b y^0}/A \; ,
\end{split}
\end{align}
 its square and that of the Cartan torsion
 \[
 I_aI^a=- \frac{16 b^2 A^{b-1}}{(b^2-1)(y^0)^{2 b} }\;, \quad C_{abc}C^{abc}=-\frac{10 b^2 A^{b-1}}{(b^2-1)(y^0)^{2 b} } .
 \]
 \end{remark}
Observe that according to our method, we have  obtained some characteristic tensors of the model pseudo-Minkowski space,
only dependent on $y^a$, and so we have just to calculate the commutation coefficients of the vierbein for those Lorentz-Finsler
spaces on which we are interested, we shall see examples in theorem \ref{afc} and in Sec.\ \ref{mvc}.

\subsection{The Berwald condition}
In this section the dimension $n+1$ is arbitrary.
Due to Eq.\ (\ref{spr}) a pseudo-Finsler space is Berwald iff
\[
H^r:=\tfrac{1}{2}c^a_{\ bc}(x) y_a(g^{rb} y^c-g^{rc}  y^b)
\]
is quadratic in $y^a$.
Given a triple $(a;b,c)$ let us say that it has property $*$ if
\[
H^{r bc}_a :=y_a(g^{rb} y^c-g^{rc}  y^b)
\]
is quadratic in $y^d$ for every $r$.
In order to prove that a space is Berwald it is sufficient to show that $c^a_{\ bc}(x) \ne 0$ only for those  triples $(a;b,c)$ for which $*$  holds true.

%%%%%%%%%

Remarkably, for the Bogoslovsky Lorentz-Minkowski space $*$ is true for every choice $(a;b,c)$ with $a>0$ and $b,c\ne n$, namely
\begin{equation}
y_a (y^0 g^{rJ}-y^Jg^{r0}), \quad y_a (y^J g^{rK}-y^Kg^{rJ})
\end{equation}
are quadratic for every $r$. This fact can be checked with a calculation but can be also understood as follows. Notice from (\ref{yaa}) that $y_a$ for $a>0$ is proportional to $l_a(y^d) F^2/A$ where $A$ is given by (\ref{aaa}) and $l_a(y^d)$ is some linear function. The inverse metric (\ref{sof}) includes a last term that cancels out when the metric is inserted into $H^{rbc}_a$, so only the other term proportional to $A/F^2$ does not cancel. The factors simplify and only a quadratic term is left.

%%%%

Less obvious is the fact that $*$ is true for $(n;n,0)$  (or equivalently $(n;0,n)$), but it can be easily checked with a calculation.
\begin{align}
\begin{split}
H_{n}{}^{0n0}=& \tfrac{(1 - b)}{1 + b}(y^0)^2\;,\quad
H_{n}{}^{Jn0} = - \tfrac{2 b}{1 + b} y^0 y^J\;,\\
H_{n}{}^{nn0}=&\tfrac{b A-(1+3b) y^0y^n}{1+b} \;.
\end{split}
\end{align}
%%%

%The fact that $H^r$ is quadratic is really due to a remarkable property of the Bogoslovsky Lorentz-Minkowski space, namely for every $k>0$ and every $r$
%\begin{equation}
%y_k (y^0 g^{r1}-y^1g^{r0}), \quad y_k (y^0 g^{r2}-y^2g^{r0}), \quad y_k (y^2 g^{r1}-y^1g^{r2})
%\end{equation}
%are quadratic in $y^a$.

%%%%%%%%%%
As a consequence, our method allows us to prove the following theorem
\begin{theorem}
%\label{fcd}
\label{afc}
Any $n+1$-dimensional Lorentz-Finsler space  of coordinates $\{u,v, \{x^i\}_{i=1,\cdots n-1}\}$ and modeled on the Bogoslovsky Lorentz-Minkowski space
\begin{equation} \label{mdx}
F=\big[2 y^0 y^n-\Sigma_{i=1}^{n-1}(y^i)^2\big]^{(1-b)/2}(y^0)^b \;,
\end{equation}
 for which the  vierbein can be chosen in the local form
\begin{align*}
y^0&=U(u)\dd u,  \\
  y^i&=\sum_{j=1}^{n-1} Z^i_j(\vec x,u) \dd x^j+ Z^i (\vec x,u) \dd u,  \quad  i=1,\ldots, n-1 \\
 y^n&= \dd v+V (\vec x,u,v) \dd u+\sum_{j=1}^{n-1} V_j(\vec x,u) \dd x^j.
\end{align*}
(where $U, Z^i_j,Z^i, V, V_j$ are arbitrary $C^1$ functions, with $Z^i_j$ invertible and $U\ne 0$) is really Berwald.
\end{theorem}

\begin{proof}
Using these vierbein an explicit computation shows that
\begin{align}
\begin{split}
&d y^0=0\;,\\
&d y^J=-\tfrac{1}{2}c^J_{KL}(x) y^K\wedge y^L-c^J_{K0}(x)y^K\wedge y^0\;,\\
&d y^n=-c^n_{J0}(x)y^J\wedge y^0-\tfrac{1}{2}c^n_{JK}(x)y^J\wedge y^K-c^n_{n0}(x)y^n\wedge y^0.
\end{split}
\end{align}
Therefore, we deduce that the property $*$ holds for the triples $(a;b,c)$ corresponding to the only non-vanishing commutation coefficients $c^a_{bc}$.
\end{proof}

% (on the contrary the expressions by Fuster and Pabst depend on both variables). At this point we have just to insert these quantities into the expressions of the notable tensors. If we were to consider a different pseudo-Finsler space based on the same Bogoslovsky pseudo-Minkowski space we  could keep  the previous calculations of the tensors $y_a,g_{ab}, g^{ab}, C_{abc}, I_a$.

\section{An example: Finsler pp-wave} \label{mvc}
We end the paper with an application of the previous theory. Recently Fuster and Pabst have proposed a Finslerian gravitational pp-wave \cite{fuster16} based on the Bogoslowsky model Lorentz-Minkowski space. Their vierbein are
%, namely $L=-\frac{1}{2} F^2$, where
%\begin{equation} \label{mdp}
%F=(2 y^0 y^3-(y^1)^2-(y^2)^2)^{(1-b)/2}(y^0)^b \;,
%\end{equation}
% $b$ is a dimensionless parameter and
\[
y^0=\dd u, \quad  y^1=\dd x, \quad  y^2=\dd y, \quad  y^3= \dd v+\tfrac{1}{2} \Phi(x,y,u) \dd u.
\]

%This is indeed a pseudo-Finsler space modeled over the same pseudo-Minkowski space, namely that introduced by  Bogoslovsky  \cite{bogoslovsky94}.
Here $\{u,x,y,v\}$ are coordinates on $M$, which would be denoted $\{x^\alpha\}$ in the notation of the previous sections, the differentials $\{\dd u, \dd x, \dd y, \dd v\}$
are nothing but the coordinates of the tangent bundle denoted by $\{y^\alpha\}$ in the previous sections,
while the previous line in display gives $\{y^a\}$, namely the vertically induced coordinates
which are used to clarify that the Finsler function is indeed based on the same Minkowski space.
In fact Eq.\ (\ref{mdp}), which has no dependence on the base variable $x$, holds true by definition.
The only non-vanishing commutation coefficients are obtained with the exterior differentiation of $y^3$, see Eq.\ (\ref{nxj}) namely
\[
c^3_{01}=-c^3_{10}=\Phi_x, \quad c^3_{02}=-c^3_{20}=\Phi_y.
\]
 At this point we have just to insert these quantities into the expressions of the notable tensors.

By our Theorem \ref{afc} the Lorentz-Finsler space is Berwald.
In fact, for the pp-wave spacetime the components of the spray computed according to equation (\ref{spr}) yield
\begin{align}
\begin{split}
G^u&=0\;,\\
G^x&=\tfrac{1}{4} (y^0)^2 \Phi_x\;,\\
G^y&=\tfrac{1}{4} (y^0)^2 \Phi_y\;,\\
G^v&=\tfrac{1}{4} y^0 \left(y^0 \Phi_u + 2 y^2\Phi_y + 2 y^1\Phi_x\right),
\end{split}
\end{align}
which are quadratic.
Since the space is Berwald, the Landsberg tensor and the Douglas curvature vanish: $L_{abc}=0$, $D^a_{bcd}=0$.

The Ricci scalar computed according to (\ref{ric}) is given by
\[
\textrm{Ric}(y)=\tfrac{1}{2} (y^0)^2 \left(\Phi_{yy} + \Phi_{xx}\right).
\]
So it vanishes iff $\Phi$ is harmonic.

%The space is Berwald because the components of $H^r=c^e_{\ m n}(x) g^{rm} y_e y^n$ entering (\ref{spr})
%are given by
%\[
%H^0=0,
%\quad H^1=\tfrac{1}{2} (y^0)^2 \Phi_x,
%\quad H^2=\tfrac{1}{2} (y^0)^2  \Phi_y,
%\quad H^3=\tfrac{1}{2} y^0 \left(y^2 \Phi_y + y^1\Phi_x\right),
%\]
%so they are quadratic in $y^a$.

%
%\section{Finslerian outgoing Eddington-Filkenstein metric}
%
%The outgoing Eddington-Filkenstein metric is
%\[
%g=
% \left(1-\tfrac{2m}{r} \right) du^2 + 2 du dr - r^2 d\Omega^2=2 y^0 y^3-(y^1)^2-(y^2)^2
%\]
%where
%\begin{align}
%y^0&=\dd u,\\
%y^1&=r\dd \theta,\\
%y^2&=r \sin \theta \dd \varphi, \\
%y^3&=\dd r+ f(r)\dd u . \qquad f(r)=\tfrac{1}{2}\left(1-\tfrac{2m}{r} \right)
%\end{align}
%We are interested in the Lorentz-Finsler metric modeled on the Boloslovsky metric with these vierbein where $f$ can be left arbitrary (possibly to be fixed imposing $Ric(y)=0$ and hence dependent on $b$).
%
%Spray? Ricci scalar?

\section{Conclusions}

We have shown that for pseudo-Finsler spaces modeled on a pseudo-Minkowski space several
Finslerian quantities can be explicit calculated. We have found interesting expressions
for the Berwald curvature, Douglas curvature, Landsberg tensor, mean Landsberg tensor,
mean Berwald curvature, and Ricci scalar in terms of: (a) vertical quantities,
which only depend on the geometry of the indicatrix, and (b) base quantities,
namely the commutation coefficients, which only depend on how the model Minkowski
space is displaced all over the manifold. These different contributions can be separately investigated
for specific metrics. As an illustration, we have applied the method to the Finslerian pp-wave recently found by Fuster and Pabst. We hope that our results could serve to bridge the gap between abstract Finsler
geometry and its applications.

\section*{Acknowledgements}
AGP wishes to thank the ``Dipartimento di Matematica e Informatica U. Dini'' at Florence University
where this work was carried out for hospitality and financial support.

AGP is supported by the projects
IT956-16 (``Eusko Jaur\-la\-ri\-tza'', Spain),
FIS2014-57956-P
(``Ministerio de Econom\'{\i}a y Competitividad'', Spain) and \\
PTDC/MAT-ANA /1275/2014 (``Funda\c{c}\~{a}o para a Ci\^{e}ncia e a Tecnologia'', Portugal). EM wishes to thank GNFM of INDAM.

%\section*{References}

%
%\bibliography{../../bibliografie/simultaneity,../../bibliografie/libri,../../bibliografie/miei,../../bibliografie/mieiPrep,../../bibliografie/mieiProc}
%\input{psefin.bbl}
%\bibliographystyle{elsarticle-num}
%\bibliographystyle{plain}
%\bibliographystyle{jmpTitles}

\end{document}